\documentclass[12pt]{amsart}
\usepackage{amsfonts,amssymb,amsthm,eucal,amsmath,verbatim,bbold,graphicx,setspace,fullpage,subfigure,lmodern}
\allowdisplaybreaks
\newcommand{\R}{\mathbb{R}}

\newcommand{\inprod}[2]{\left\langle #1, #2 \right\rangle}

\renewcommand{\P}{\mathbb{P}}

\begin{document}

\title{Comparison of $f$-vectors of Random Polytopes to the
  Gaussian Distribution}
\author{Sang Du and Mark Syvuk}

\maketitle

\begin{abstract}
Choose $n$ random, independent points in $\R^d$ according to a fixed
distribution. The convex hull of these points is a random
polytope. In some cases, central limit theorems have been proven
for the components of $f$-vectors of random polytopes constructed in
this and similar ways. In this paper, we provide numerical evidence
that the components of the $f$-vectors of random
polytopes generated according to five different distributions are
approximately jointly Gaussian for large $n$.
\end{abstract}


\section{Introduction}
Phenomena in high dimensions have been the subject of many
studies. In particular, one may be interested in studying higher
dimensional geometric objects called polytopes. A $polytope$ $P \in
\R^d$ is a bounded subset of $\R^d$ formed by taking the intersection
of finitely many halfspaces. The $f$-vector of a polytope $P$ is a
vector $f = (f_0, f_1 . . ., f_{d-1})$ such that $f_i$ is the number of $i$-dimensional
faces of $P$.  In particular, $f_0$ is the number of vertices of $P$
and $f_{d-1}$ is the
number of $(d-1)$-dimensional faces, known as $facets$. \par There are
some known results on the distribution of the $f$-vector of types of
random polytopes. Previous work has
found that individual components of the $f$-vectors are approximately
Gaussian when
the polytopes are generated according to various distributions. In
particular, B\'ar\'any and Vu \cite{clt_for_gaussian} proved that if a polytope is taken as the
convex hull of independent and identically distributed (i.i.d.) points
according to the standard Gaussian distribution, then the volume and
individual components of the $f$-vector of that polytope satisfy a
central limit theorem; that is, converge to a Gaussian distribution
as the number of i.i.d. random points tends to infinity. B\'ar\'any and
Reitzner \cite{poisson} proved
central limit theorems for the volume and components of the
$f$-vector of a so-called Poisson random polytope; that is, the convex hull of
the intersection of an arbitrary volume 1 convex body $K$ with a
Poisson process $\chi$. Following earlier work of Reitzner \cite{reitz_shit}, Vu \cite{clt_for_rand_poly} proved that the volume and
components of the $f$-vector of a random polytope drawn from points
inside a smooth convex set also satisfy central limit
theorems. \par In this paper, we consider the joint
distributions of the $f$-vectors of random polytopes. The method used
in this paper for constructing a random polytope is taking the convex
hull of a set of i.i.d. random points generated from a fixed
distribution. Our conjecture is the following:

\subsection{Conjecture} \emph{Let N be a positive integer. Let
  $\mathcal{X}_N = \{X_1, . . ., X_n\}$ be a collection of
  random points $X_i
  \in \R^d$, which are independent and identically distributed according to a
  fixed distribution.  Letting $P_N \in \R^d$ be the convex hull
  of $\mathcal{X}_N$, under mild assumptions on the distribution of the
    $X_i$, the joint distribution of the $f$-vector satisifes a
    central limit theorem}.

\subsection{Results}
In the five underlying distributions of the random points
considered, numerical evidence supports the conjecture. The metric
used to compare the closeness of the distribution of the
$f$-vectors with the Gaussian distribution is based on the Kolmogorov
distance, $D_K$. For now one may think of $D_K$ as a tool that
measures the distance between two distributions and that $0 \le D_K \le 1$. If $D_K$ is a small value
then the two distributions are considered to be close to each
other. The values of $D_K$ presented in the table are for comparison
to the appropriate Gaussian distribution; details are discussed in
section 3. Here, $d$ is the dimension, $n$ is the number of random points, and $N$ is the sample size or number of $f$-vectors.
\\ \par 

\begin{center}
\begin{table}[ht]
\caption{Uniform Distribution in Cube}
\begin{tabular}{ccc|r@{.}l}
$d$ & $n$   & $N$   & \multicolumn{2}{c}{$D_K$} \\
\hline
5   & 64000 & 25000 & 0 & 008995 \\
6   & 64000 & 25000 & 0 & 007590 \\
7   & 8000  & 25000 & 0 & 007133 \\
8   & 1000  & 3125  & 0 & 02613
\end{tabular}
\label{table:cube}
\end{table}

\begin{table}[ht]
\caption{Uniform Distribution in $\ell_1$ Ball}
\begin{tabular}{ccc|r@{.}l}
$d$ & $n$   & $N$   & \multicolumn{2}{c}{$D_K$} \\
\hline
5   & 64000 & 25000 & 0 & 009971 \\
6   & 64000 & 25000 & 0 & 01075 \\
7   & 4000  & 25000 & 0 & 01113 \\
8   & 1000  & 2000  & 0 & 03003
\end{tabular}
\label{table:l1}
\end{table}

\begin{table}[ht]
\caption{Uniform Distribution in $\ell_2$ Ball}
\begin{tabular}{ccc|r@{.}l}
$d$ & $n$   & $N$   & \multicolumn{2}{c}{$D_K$} \\
\hline
5   & 64000 & 25000 & 0 & 008495 \\
6   & 64000 & 25000 & 0 & 007873 \\
7   & 2000  & 25000 & 0 & 01182 \\
8   & 1000  & 940   & 0 & 04573
\end{tabular}
\label{table:l2}
\end{table}

\begin{table}[ht]
\caption{Standard Gaussian Distribution}
\begin{tabular}{ccc|r@{.}l}
$d$ & $n$   & $N$   & \multicolumn{2}{c}{$D_K$} \\
\hline
5   & 64000 & 25000 & 0 & 01197 \\
6   & 64000 & 25000 & 0 & 01125 \\
7   & 64000 & 25000 & 0 & 009252 \\
8   & 32000 & 1000  & 0 & 01530
\end{tabular}
\label{table:gaussian}
\end{table}

\begin{table}[ht]
\caption{Uniform Distribution in Hemisphere}
\begin{tabular}{ccc|r@{.}l}
$d$ & $n$   & $N$   & \multicolumn{2}{c}{$D_K$} \\
\hline
5   & 64000 & 25000 & 0 & 006166 \\
6   & 4000 & 25000 & 0 & 008827 \\
7   & 500 & 12500 & 0 & 02883 \\
\end{tabular}
\label{table:gaussian}
\end{table}

\end{center}
Since $D_K$ is a metric that takes on values in
the interval [0,1], a value with order of magnitude $10^{-3}$ as the Kolmogorov distance
should be considered a rather small value. \par
Figure 1 contains histograms of individual components of 64,000 $f$-vectors computed from the
convex hull of 25,000 i.i.d. random points uniformly distributed in the
5-dimensional cube. After standardizing the
$f$-vectors to have sample mean zero and identity sample covariance,
the resulting data lie in $\R^2$. As further illustration, the 2-dimensional data of the standardized
$f$-vector is plotted in $\R^3$. 
\begin{figure}[htp]
\subfigure[$f_0$]{\includegraphics[width=1.5in]{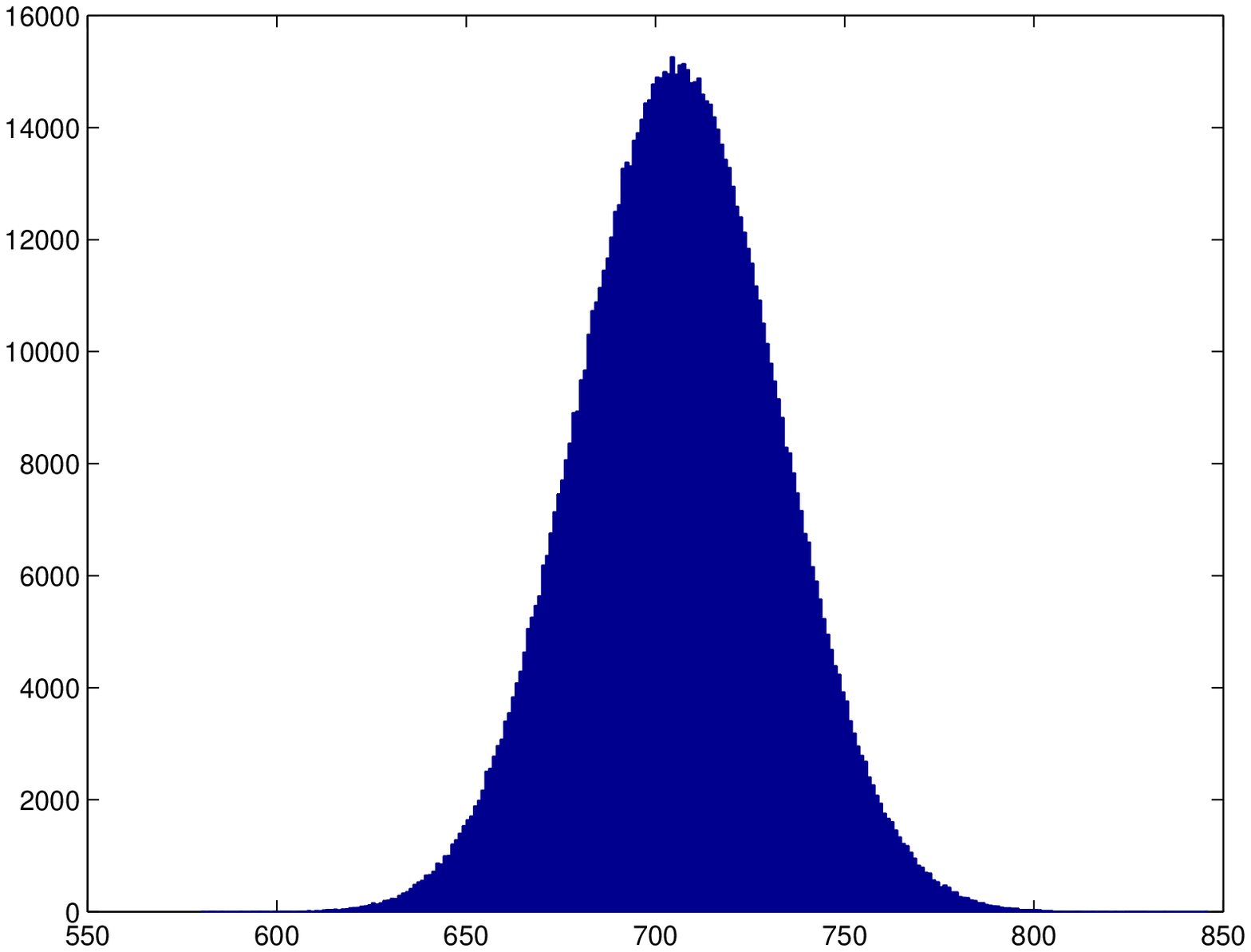}}
\subfigure[$f_1$]{\includegraphics[width=1.5in]{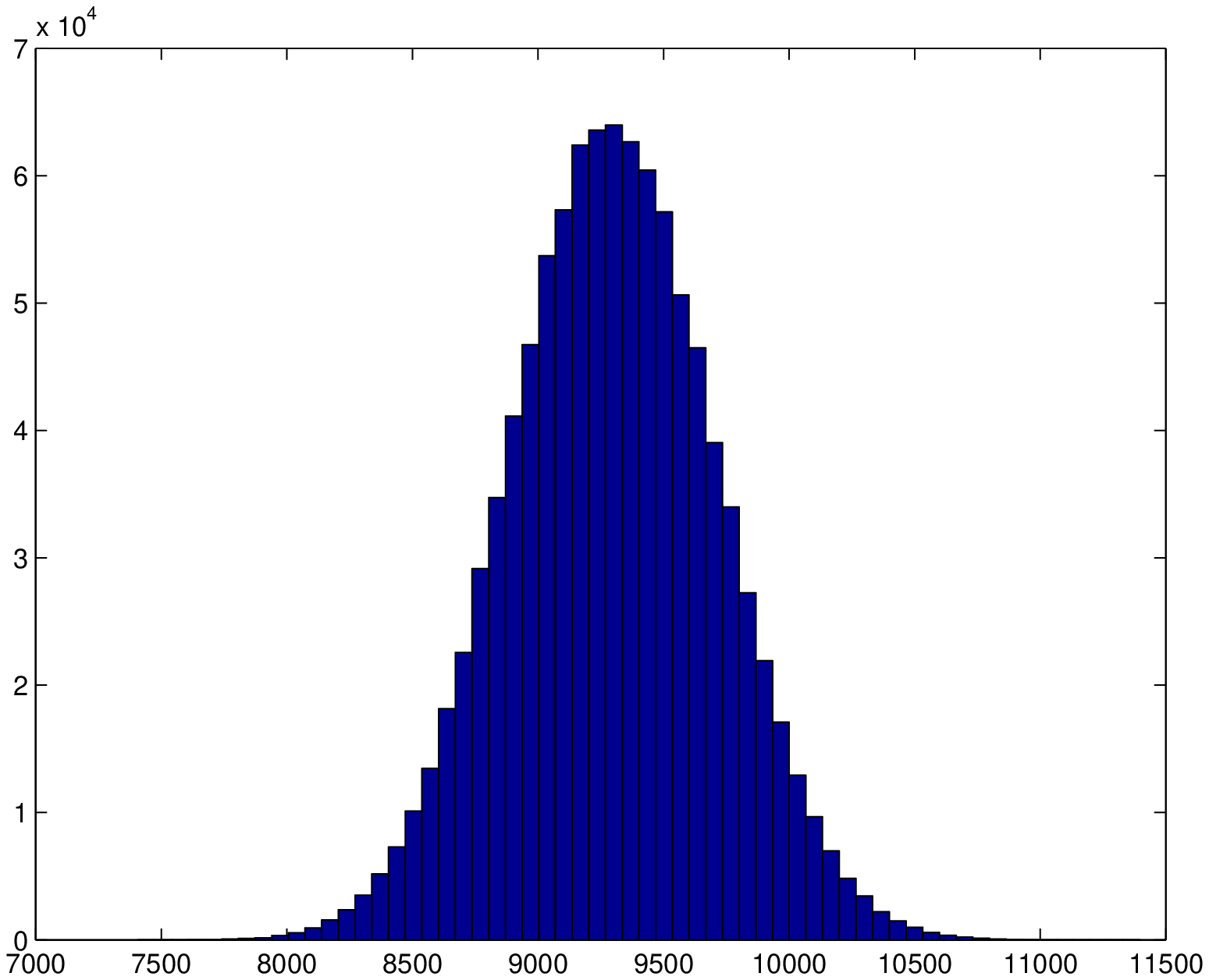}}
\subfigure[$f_2$]{\includegraphics[width=1.5in]{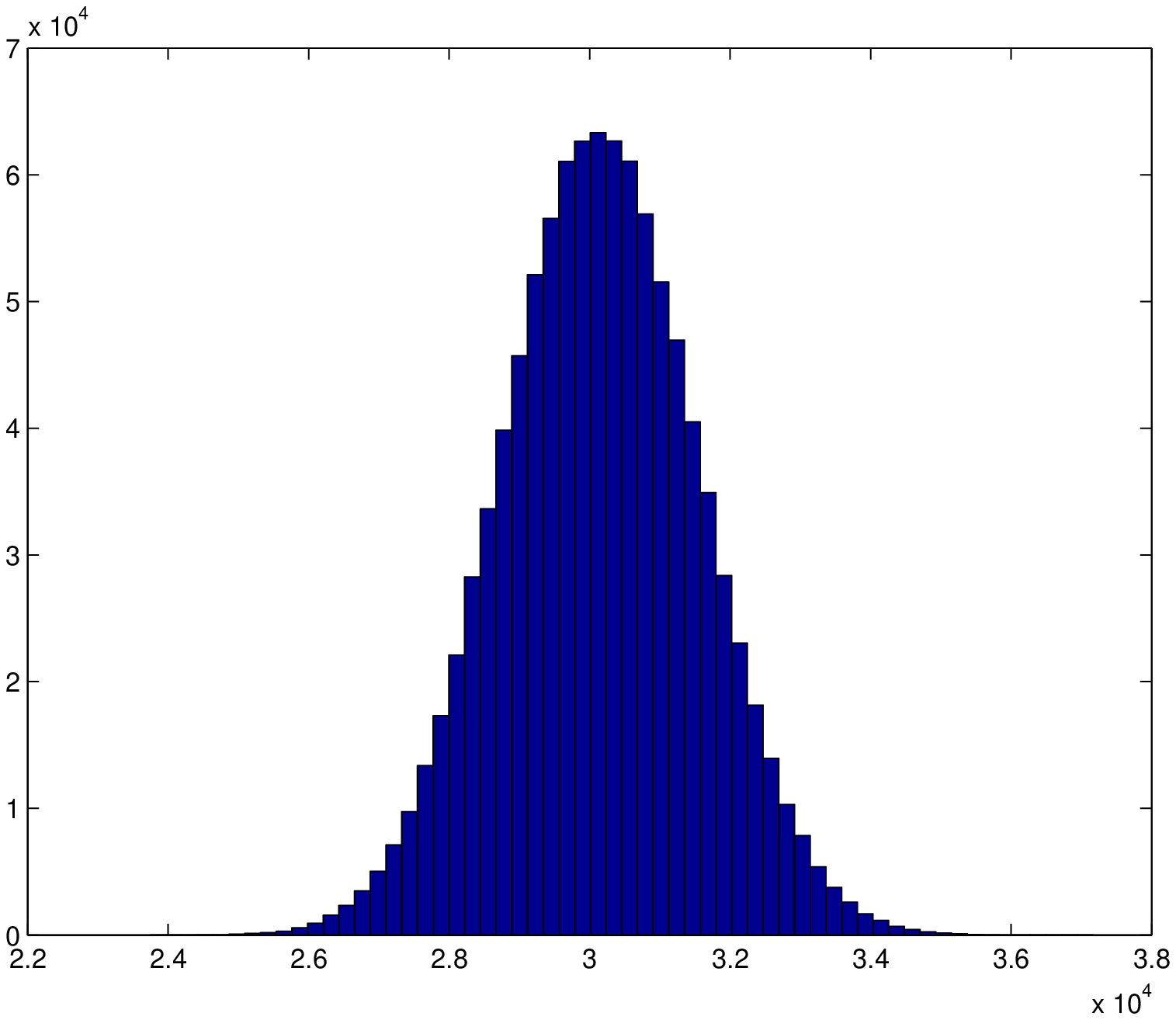}}
\subfigure[$f_3$]{\includegraphics[width=1.5in]{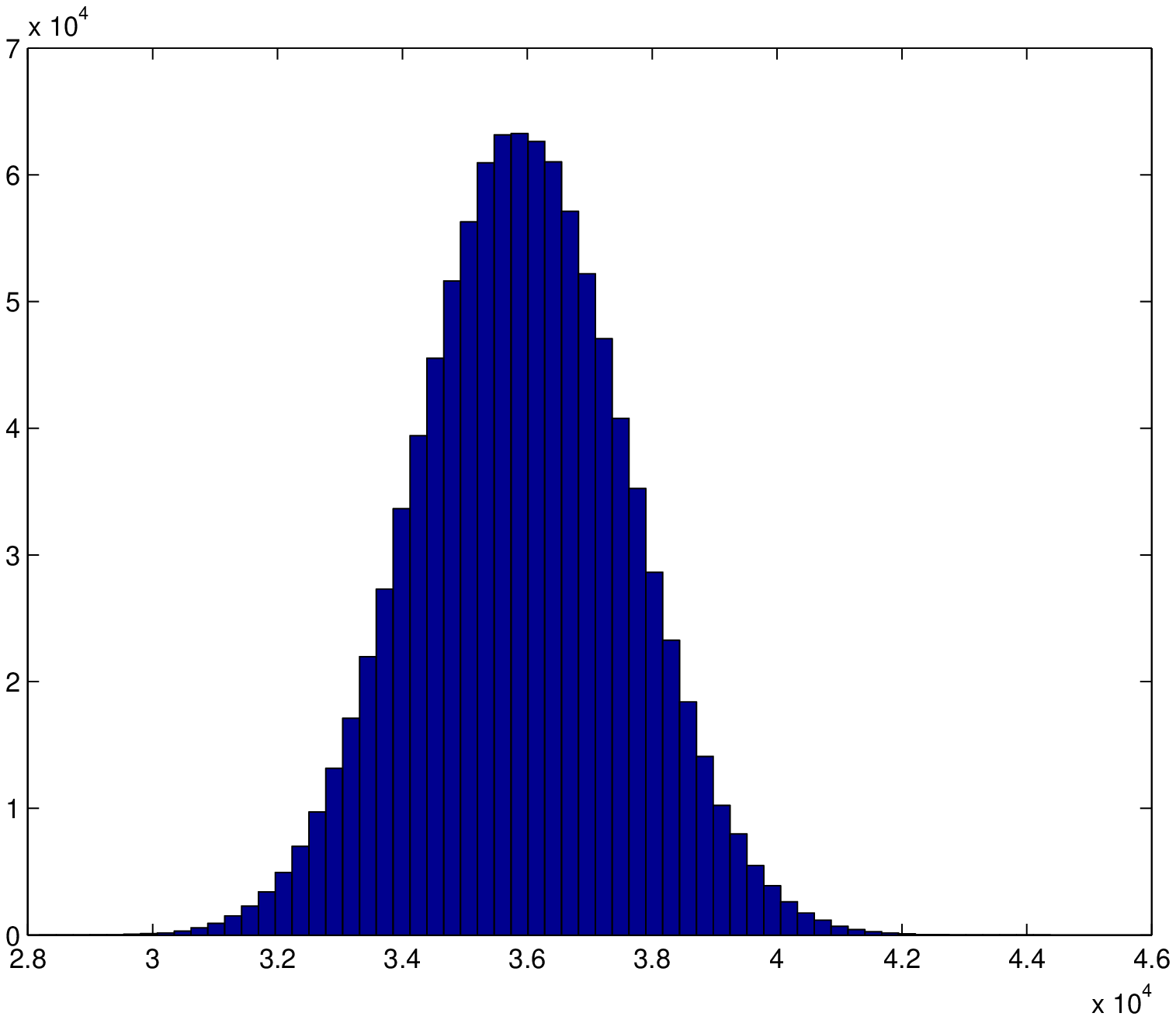}} 
\subfigure[$f_4$]{\includegraphics[width=1.5in]{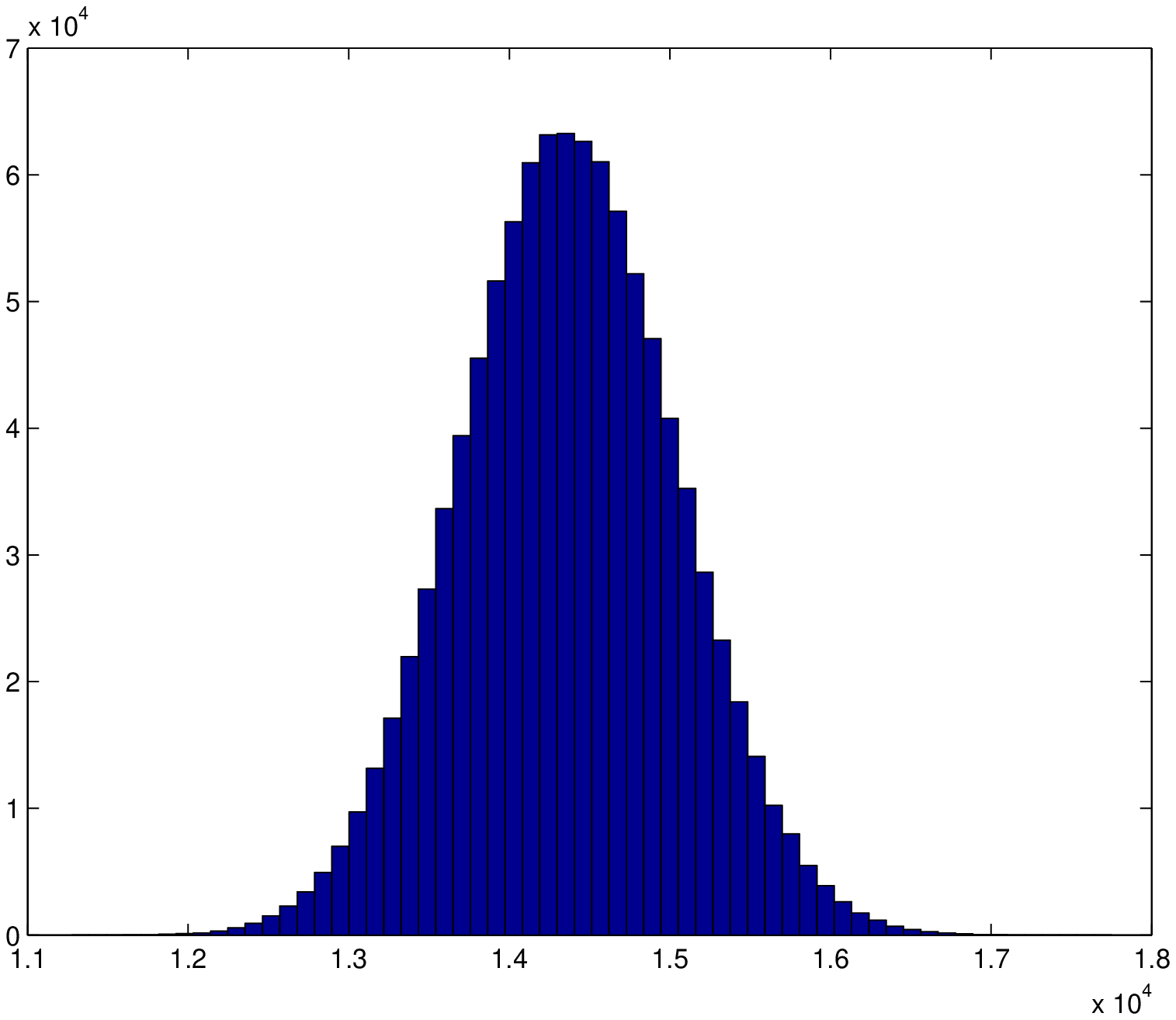}} 
\caption{Histograms of components of $f$-vectors in $\R^5$. Data came
  from 25,000 $f$-vectors computed from random polytopes generated by
  64,000 i.i.d. uniformly distributed points in the 5-dimensional cube.}
\end{figure}


\begin{figure}
\begin{center}
\includegraphics[scale=.6]{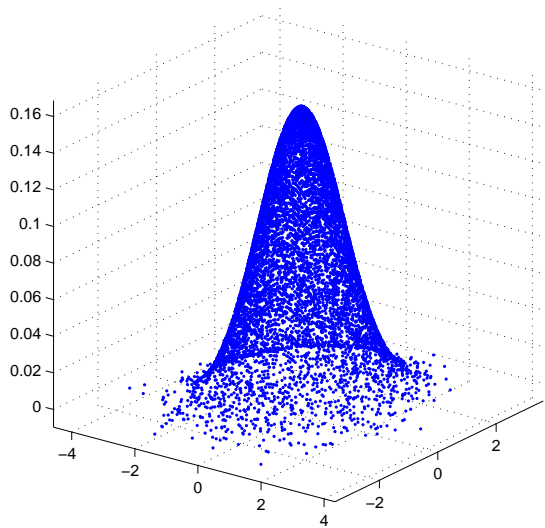} 
\end{center}  
\caption{Standardized $f$-vectors
  in $\R^5$. Data came
  from 25,000 $f$-vectors computed from random polytopes generated by
  64,000 i.i.d. uniformly distributed points in the 5-dimensional cube.} 
\end{figure}

\section{Background}
It was mentioned in the introduction that one way to construct a random polytope is to
take the convex hull of random points. One way to generate random
points is to let
$\{X_i\}_{i=1}^n$ be a collection of points generated independently
according to a fixed distribution. Five distributions used to generate
random points were considered and the random polytope
studied was the convex hull of the $\{X_i\}_{i=1}^n$. The underlying
distributions considered are the following:\par 
\begin{enumerate}
\item Uniform in the Cube. Define
$B_{\infty}^d := \{x \in \R^d: ||x||_\infty \le 1 \}$  where
$||x||_\infty := \underset{1 \le i \le d}{max} |x_i|$. A random point $X \in B_{\infty}^d$ is
constructed by picking each component of $X$ independently according to the uniform
distribution in the interval [0,1]. The resulting random point $X$ is
uniformly distributed in $B_{\infty}^d$. 
\item Uniform in the Euclidean Ball. Define $B_2^d$
in the $\ell_2$ norm as
$B_2^d:= \{x \in \R^d: ||x||_2 \le 1\}$ where $||x||_2:= \left( \sum\limits_{i=1}^d
  |x_i|^2 \right) ^{\frac{1}{2} }$. To generate a random point $X \in
\mathbb{R}^d$ such that $X$ is uniformly distributed in $B_2^d$, first
generate a point $Z = (Z_1,. . ., Z_{d+2})$ that is distributed
according to the standard $d+2$ dimensional Gaussian
distribution. This is done by picking each component of $Z$ independently
according to the standard univariate normal distribution.  Let
$X = \frac{(Z_1, . . ., Z_d)}{||Z||_2}$, where $||Z||_2 = \left(
  \sum\limits_{i=1}^{d+2} |Z_i|^2 \right) ^{\frac{1}{2} }$. Then it is
classical that $X$ is uniformly distributed in
$B_2^d$. See [1] for a proof.
\item Uniform in the $\ell_1$ ball.
 Define $B_1^d$ in the $\ell_1$ norm as $B_1^d:= \{x \in
  \R^d: ||x||_1 \le 1\}$ where  $||x||_1:=
\sum\limits_{i=1}^d |x_i|$. Let $Z= (Z_1, . . ., Z_{d+1})$ where the
components of $Z$ are generated independently according to the
distribution with the density
function $f(x)=\frac{1}{2}e^{-|x|}$.  Let
$X = \frac{(Z_1, . . ., Z_d)}{||Z||_1}$, where $||Z||_1 = 
  \sum\limits_{i=1}^{d+1} |Z_i| $; it is
proven in \cite{balls} that $X$ is uniformly distributed in
$B_1^d$. 
\item The Standard Normal distribution.
Let $X \in \R^d$ be a random vector and generate each components of
$X$ independently
according to the standard normal distribution. Then $X \in
\mathbb{R}^d$ is distributed as a standard normal random
vector. 
\item Uniform in the Hemisphere. Similar to (2), but take the absolute value of
  the first component of the random point $X$ in (2). Note that unlike
  the underlying convex bodies in (1)-(3), this body is neither
  smooth nor a polytope.
\end{enumerate}
\par


\section{Computation and Analysis of \emph{f}-vectors}
The software MATLAB and QHULL were used to carry out the following simulations and
computations. 
One realization of
the $f$-vector is computed by first fixing the dimension $d$, generating
$n$ i.i.d. points  $\{X^i\}_{i=1}^n$ from one of the fixed probability
distributions described above,
taking the convex hull of these points to construct a random polytope $P$,
and computing the $f$-vector of the random polytope.
Assuming the random polytope is \emph{simplicial}, that is, each of its facets has exactly $d$ vertices, the $f$-vector
of the polytope can be directly computed from the facets. Let $e_1, ... e_m \subset \R^d$ be the $m$ facets of $P$,
where $e_{i,j}$ is the $j$th vertex of the $i$th facet, for $1\leq i\leq m$,
$1\leq j\leq d$. By the definition of $f$, $f_{d-1}=m$. Since $f_{d-2}$ is the number of $(d-2)$-dimensional
faces of $P$, $f_{d-2}$ is found by counting all intersections of size $(d-1)$ between pairs of
facets $e_i$ and $e_j$, for $i\neq j$. Since each $(d-2)$-dimensional face is the intersection of exactly two
facets, this accounts for all $(d-2)$-dimensional faces. For example, 3-dimensional faces $(x_1,x_2,x_3)$ and $(x_1,x_2,x_4)$,
where each $x_i$ is a vertex for that face,
have intersection $(x_1,x_2)$ which is of size 2 so it must be an edge, but $(x_1,x_2,x_3)$ and $(x_1,x_4,x_5)$ have
intersection $(x_1)$ which is of size 1 so it is not an edge. Let $g_1, ..., g_k\in\R^{d-1}$ be these distinct
intersections of size $(d-2)$; clearly $f_{d-2} = k$. By keeping track of all the intersections
of size $(d-1)$, we obtain a list of all $(d-2)$-dimensional faces of $P$. We can then repeat
this process with those faces. Continuing inductively in this manner, the remaining
components of $f$ can be determined.

Note that the convex hull of points in general can be non-simplicial
(a facet may have more than $d$ vertices), and so the above algorithm
for constructing the $f$-vector from the facets is not valid for all
polytopes. However, the probability of a non-simplicial convex hull
arising from the distributions of points considered in this paper is
0. \par
A collection of $N$ realizations of $f$-vectors was obtained, and the
data were then standardized and compared to the standard Gaussian
distribution. The following demonstrates the computations of
standardizing the $f$-vectors. Let $f^k$
be the $k^{th}$ $f$-vector; then $f^k_i$ is the $i^{th}$ component of
the $k^{th}$ $f$-vector. Define the sample mean $\bar{f}$\ of the
\emph{f}-vectors as $\bar{f}= (
\bar{f}_0, \bar{f}_1, . . ., \bar{f}_{d-1 } )$ where
$\bar{f}_j = \frac{1}{N}  \sum\limits_{k=1}^N f^k_j
$. Given
the sample data $\{f^k\}_{k=1}^N$, the sample covariance matrix $S=[s_{ij}]_{i,j=1}^d$ is
a $d$-by-$d$ matrix with entries given by
$s_{ij} = \frac{1}{N-1}\sum\limits_{k=1}^N (f^k_i -
\bar{f_i})(f^k_j - \bar{f_j})$. To standardize to identity sample
covariance, first
diagonalize the sample covariance $S= UDU^T$, where $D$ is a diagonal matrix of
eigenvalues and $U$ is an orthogonal matrix where the columns of $U$ are
the eigenvectors. Moreover, the eigenvalues of $D$ can be assumed to be in decreasing
order, $\lambda_1 \ge \lambda_2 \ge . . . \lambda_d$.
Note that the
components of the $f$-vectors have some dependence on each
other. For polygons, the components of the $f$-vector obey $f_0 = f_1$. For
convex polyhedra in $\R^3$, the components of the $f$-vector satisfy $f_0 - f_1
+ f_2 = 2$, known as Euler's relation. In higher dimensions, linear
dependence
between the components of the $f$-vector causes the covariance matrix
of the $f$-vector to be singular. These type of linear dependencies are
known as Dehn-Sommerville equations for simplicial polytopes [6]. A singular covariance matrix leads to zero eigenvalues in the diagonalization
of the covariance matrix. In practice, the distinction between
non-zero eigenvalues and approximately zero eigenvalues of the sample
covariance matrix was easy to
identify. In such case, let the eigenvalues of the
sample covariance matrix that are 
close to zero be eliminated to ensure that $D^{-\frac{1}{2}}$ is
well-defined. So define a new matrix $D^*  = \begin{bmatrix}
 \lambda_1 & \cdots & 0 \\
\vdots &\ddots &   \vdots \\
0 & \cdots & \lambda_p \\
\end{bmatrix} $
where $D^*$ is a $p$-by-$p$ matrix with the
eigenvalues of matrix $D$ that are close to zero eliminated. Also
define a new matrix $U^* = \begin{bmatrix}  
| &  & |  \\ 
u_1  & \cdots & u_p  \\ 
| & & | \\ 
\end{bmatrix}$ where the eigenvector $u_i$ corresponds to the $\lambda_i$ eigenvalue.
Then let $\hat{f^k} = (D^*)^{-\frac{1}{2}}(U^*)^T \tilde{f^k}$, where $\tilde{f^k}$ is the $k^{th}$
re-centered $f$-vector. Then $\{\hat{f}^k\}_{k=1}^N \subseteq \R^p$ is now a
standardized data set with mean 0 and identity sample covariance. \par

After standardization, one can now compare the standardized $f$-vectors to the standard Gaussian distribution, in terms of the Kolmogorov distance.  The Kolmogorov distance between two random variables $X,Y$ is defined to be $d_K:=\sup_{t\in\R}|\P[X\le t]-\P[Y\le t]|$.  Observe that since the standardized data are conjectured to be approximately Gaussian, any linear combination of their components should also be close to Gaussian.  This motivates the following numerical test.

Fix $M$, and let $\{a_i\}_{i=1}^M \subseteq \mathbb{S}^{p-1}$ be $M$ independent samples of uniformly chosen points on the sphere in $\R^p$. For each $a_i$, compute the Kolmogorov distance
\[d_K^i:=\sup_{t\in\R}|F_{N,i}(t)-\Phi(t)|\]
where $F_{N,i}(t)$ is the empirical cumulative distribution function of the data projected onto the direction $a_i$; that is,
\[F_{N,i}(t):=\frac{1}{N}\#\{\inprod{a_i}{\hat{f}^k}\le t\},\] and $\Phi(t)$ is the cumulative distribution function of the standard normal distribution.  Finally, define a measure $D_K$ of the distance from the data to Gaussian by
\[D_K:=\sup_{1\le i\le M}d_K^i.\] The values presented in the tables in Section 1 are obtained using
$M=100,000$.

\subsection{Acknowledgements}
The authors would like to thank advisor Elizabeth Meckes for her
guidance and mentoring on the subject along with her help in editing this paper. This work was supported by NSF
grants DMS-0905776 and DMS-0852898.

\end{document}